\documentclass[11pt,twoside]{article}
\usepackage{mathbbold}
\usepackage{amsmath, amssymb, mathrsfs, graphicx}

\def\scr{\mathscr}

\def\az{\alpha}  \def\bz{\beta}
    
\def\ez{\eta}    \def\fz{\varphi}
\def\gz{\gamma}  \def\kz{\kappa}
\def\lz{\lambda} 
     \def\oz{\omega}

        \def\sz{\sigma}
        \def\uz{\theta}
\def\vz{\varepsilon} \def\xz{\xi}

\def\llz{\Lambda}

\def\qd{\quad}
\def\qqd{\qquad}

\newcommand{\mathsym}[1]{{}}

\def\scr{\mathscr}

\def\le{\leqslant}
\def\ge{\geqslant}

\font\cms=cmss9 scaled \magstep1

\def\nnd{\noindent}

\def\thm{\nnd\bg{thm1}}
\def\crl{\nnd\bg{crl1}}

\def\prp{\nnd\bg{prp1}}

\def\xmp{\nnd\bg{xmp1}}

\def\dethm{\end{thm1}}
\def\decrl{\end{crl1}}

\def\deprp{\end{prp1}}
\def\dexmp{\end{xmp1}}

\def\prf{\medskip \noindent {\bf Proof}. }
\def\deprf{\quad $\square$ \medskip}
\def\bg{\begin}
\def\be{\bg{equation}}
\def\de{\end{equation}}

\def\dear{\end{eqnarray}}
\def\lb{\label}
\def\ct{\cite}

\newcommand{\rf}[2]{[\ref{#1}; #2]}

\def\den{\end{enumerate}}

\DeclareMathOperator*{\var}{Var} 
\def\d{\text{\rm d}}

\begin{document}

\allowdisplaybreaks[4]
\thispagestyle{empty}
\renewcommand{\thefootnote}{\fnsymbol{footnote}}

\noindent {Chapter 6 in ``Probability Approximations and Beyond'', Lecture Notes in Statistics 205, 2012}

\vspace*{.5in}
\begin{center}
{\bf\Large Basic Estimates of Stability Rate for One-dimensional Diffusions}
\vskip.15in {Mu-Fa Chen}
\end{center}
\begin{center} (Beijing Normal University, Beijing 100875, China) \end{center}
\vskip.1in

\markboth{\sc Mu-Fa Chen}{\sc Basic estimates of exponential rate}


\date{}


\footnotetext{{\it Key words and phases}.\quad
First nontrivial eigenvalue, Hardy's inequality, one-dimensional diffusion,
principal eigenvalue, Poincar\'e-type inequality, stability rate.}
\footnotetext{In the published version as Chapter 6 in the book, Theorem 2.1, Proposition 3.2,
Corollary 4.3 and Example 5.1 here are relabeled as Theorem 6.1, Proposition 6.1,
Corollary 6.1 and Example 6.1, respectively. Similarly, formulas (1)--(36) here are relabeled as (6.1)--(6.36).
}
\bigskip

\begin{abstract}
In the context of one-dimensional diffusions, we present
basic estimates (having the same lower and upper bounds with a factor of 4 only)
for four Poincar\'e-type (or Hardy-type) inequalities. The derivation of two estimates
have been open problems for quite some time. The bounds provide exponentially ergodic or decay rates. We refine the bounds and illustrate them with typical examples.
\end{abstract}

\medskip

\section{Introduction}
An earlier topic on which Louis Chen has studied is about the
Poincar\'e-type inequalities (see \ct{clh85, clh87}, for
instance). We now use this good chance to introduction in Section \ref{s0} some recent
progress on the topic, especially for one-dimen\-sional diffusions
(elliptic operators). The basic estimates of exponentially ergodic
(or decay) rate and the principal eigenvalue in different cases are
presented. Here the term ``basic'' means that upper and lower bounds are
given by an isoperimetric constant up to a factor four. As a
consequence, the criteria for the positivity of the rate and the
eigenvalue are obtained. The proof of the main result is sketched in Section \ref{s2}.
The materials given in Sections \ref{s3}, \ref{s4}, and
Appendix are new. In particular, the basic estimates are
refined in Section \ref{s3} and the results are illustrated through examples
in Section \ref{s4}. The coincidence of the exponentially decay rate and
the corresponding principal eigenvalue is proven in Appendix for a
large class of symmetric Markov processes.

\section{The main result and motivation}\lb{s0}
\subsection{Two types of exponential convergence}

Let us recall two types of exponential convergence often studied for
Markov processes. Let $P_t(x, \cdot)$ be a transition probability on
a measurable state space $(E, {\scr E})$ with stationary
distribution $\pi$. Then the process is called {\it exponentially
ergodic} if there exists a constant $\vz> 0$ and a function $c(x)$
such that \be \|P_t(x, \cdot)-\pi\|_{\var}\le c(x) e^{-\vz t},\qqd
t\ge 0,\; x\in E.\lb{1}\de Denote by $\vz_{\max}$ be the maximal
rate $\vz$. For convenience, in what follows, we allow
$\vz_{\max}=0$. Next, let $L^2(\pi)$ be the real $L^2(\pi)$-space
with inner product $(\cdot, \cdot)$ and norm $\| \cdot \|$
respectively, and denote by $\{P_t\}_{t\ge 0}$ the semigroup of the
process. Then the process is called to have {\it $L^2$-exponential
convergence} if there exists some $\ez\,(\ge 0)$ such that \be \|P_t
f-\pi(f)\|\le \|f-\pi(f)\| e^{-\ez t}, \qqd t\ge 0,\; f\in
L^2(\pi),\lb{2}\de where $ \pi(f)=\int_E f\d \pi$. It is known that
$\ez_{\max}$ is described by $\lz_1$: \be \lambda_1=\inf\{(f, - L
f): f\in {\scr D}(L),\; \pi(f)=0,\; \|f\|=1\},\de where $L$ is the
generator with domain ${\scr D}(L)$ of the semigroup in $L^2(\pi)$.
Even though the topologies for these two types of exponential
convergence are rather different, but we do have the following
result.

\thm[\ct{cmf91, cmf05a}]{\cms For a reversible Markov process
with symmetric measure $\pi$, if with respect to $\pi$, the
transition probability has a density $p_t(x, y)$ having the property
that the diagonal elements $p_s(\cdot, \cdot)\in L^{1/2}_{\rm
loc}(\pi)$ for some $s>0$, and a set of bounded functions with
compact support is dense in $L^2(\pi)$, then we have {$\vz_{\max}
=\lz_1$}}. \dethm

As an immediate consequence of the theorem, we obtain some criterion
for $\lz_1>0$ in terms of the known criterion for $\vz_{\max}>0$.
In our recent study, we go to the opposite direction:
estimating $\vz_{\max}$ in terms of the spectral theory.

We are also going to handle with the non-ergodic case in which (\ref{2}) becomes
\be \mu\big((P_t f)^2\big)\le \mu\big(f^2\big) e^{-2\ez t},
\qqd t\ge 0,\; f\in L^2(\mu),\de
where $\mu$ is the invariant measure of the process.
Then $\ez_{\max}$ becomes
\be \lambda_0
=\inf\big\{-\mu (f L f): f\in {\scr C},\;  \mu\big(f^2\big)=1\big\},\de
where ${\scr C}$ is a suitable core of the generator, the smooth functions
with compact support for instance in the context of diffusions.
However, the totally variational norm in (\ref{1}) may be meaningless
unless the process being explosive.
Instead of (\ref{1}), we consider the following exponential convergence:
\be P_t(x, K) \le c(x, K)\, e^{-\vz t},\qqd t\ge 0,\; x\in E, \; K\!: \mbox{compact},\de
where for each compact $K$, $c(\cdot, K)$ is locally $\mu$-integrable.
Under some mild condition,
we still have $\vz_{\max}=\lz_0$. See Appendix for more details.

\subsection{Statement of the result}

We now turn to our main object: one-dimensional diffusions. The state
space is $E:=(-M, N)$\, $(M, N\le \infty)$. Consider an elliptic operator
$$L=a(x)\frac{\d^2}{\d x^2}+ b(x)\frac{\d}{\d x},$$
where $a>0$ on $E$. Then define a function $C(x)$ as follows:
$$C(x)=\!\int_{\uz}^x \frac b a,\qqd x\in E,$$
where $\uz\in E$ is a reference point. Here and in what follows, the
Lebesgue measure $\d x$ is often omitted. It is convenient for us to
define two measures $\mu$ and $\nu$ as follows.
$$ \mu(\d x)=\frac{e^{C(x)}}{a(x)}\d x,\qqd
\nu(\d x)= e^{-C(x)}\d x.
$$
The first one has different names: {\it speed}, or {\it invariant},
or {\it symmetrizable measure}. The second one is called {\it scale
measure}. Note that $\nu$ is infinite iff the process is recurrent.
By using these measures, the operator $L$ takes a very compact form
\be L=\frac{\d}{\d \mu}\, \frac{\d}{\d \nu}\qqd\Big(\text{i.e.,}\qd
Lf\equiv a\, e^{-C}\big(f' e^{C}\big)'\Big)  \lb{7}\de
which goes back to a series of papers by W. Feller, for instance \ct{feller}.

Consider first the special case that $M, N< \infty$. Then the ergodic case means
that the process has reflection boundaries at $-M$ and $N$. In analytic language,
we have Neumann boundaries at $-M$ and $N$: the eigenfunction $g$ of
$\lz_1$ satisfies $g'(-M)=g'(N)=0$.
Otherwise, in the non-ergodic case, one of the boundaries becomes
absorbing. In analytic language,
we have Dirichlet boundary at $-M$ (say): the eigenfunction $g$ of
$\lz_0$ satisfies $g(-M)=0$. Let us use codes ``D'' and ``N'', respectively,
to denote the Dirichlet and Neumann boundaries. The corresponding
minimal eigenvalues of $-L$ are listed as follows.
\begin{itemize}   \setlength{\itemsep}{-0.8ex}
\item $\lz^{\text{\rm NN}}$: Neumann boundaries at $-M$ and $N$,
\item $\lz^{\text{\rm DD}}$: Dirichlet boundaries at $-M$ and $N$,
\item $\lz^{\text{\rm DN}}$: Dirichlet at $0$ and Neumann at $N$,
\item $\lz^{\text{\rm ND}}$: Neumann at $0$ and Dirichlet at $N$.
\end{itemize}
We call them {\it the first non-trivial} or {\it the principal
eigenvalue}. In the last two cases, setting $M=0$ is for convenience
in comparison with other results to be discussed later but it is not
necessary. Certainly, this classification is still meaningful if $M$
or $N$ is infinite. For instance, in the ergodic case, the process
will certainly come back from any starting point and so one may
imagine the boundaries $\pm \infty$ as reflecting. In other words,
the probabilistic interpretation remains the same when $M$,
$N=\infty$. However, the analytic Neumann condition that
$\lim_{x\to\pm\infty} g'(x)=0$ for the eigenfunction $g$ of
$\lz^{\text{\rm NN}}$ may be lost (cf. the first example given in
Section \ref{s4}). More seriously, the spectrum of the operator may
be continuous for unbounded intervals. This is the reason why we
need the $L^2$-spectral theory. In the Dirichlet case, the analytic
condition that $\lim_{x\to\pm\infty} g(x)=0$ can be implied by the
definition given below, once the process goes $\pm\infty$
exponentially fast. Now, for general $M, N\le \infty$, let
$$\aligned
D(f)&=\int_{-M}^N {f'}^2 e^C,\qqd M, N\le \infty,\; f\in {\scr A}(-M, N),\\
 {\scr A}(-M, N)&= \text{the
set of absolutely continuous functions
on $(-M, N)$},\\
{\scr A}_0(-M, N)&= \{f\in {\scr A}(-M, N): f\text{ has a compact
support} \}.
\endaligned
$$
From now on, the inner product
$(\cdot, \cdot)$ and the norm $\|\cdot\|$ are taken with respect to
$\mu$ (instead of $\pi$).
Then the principal eigenvalues are defined as follows.
\begin{align}
\lz^{\text{\rm DD}}&=\inf\{D(f): f\in {\scr A}_0(-M, N), \; \|f\|=1\},\\
\lz^{\text{\rm ND}}&=\inf\{D(f): f\in {\scr A}_0(0, N), \; f(N-)=0\text{ if }N<\infty,\; \|f\|=1\},\\
\lz^{\text{\rm NN}}&=\inf\{D(f): f\in {\scr A}(-M, N), \;\mu(f)=0,\;
\|f\|=1\},\\
\lz^{\text{\rm DN}}&=\inf\{D(f): f\in {\scr A}(0, N), \; f(0+)=0,\;\|f\|=1\}.\lb{11}
\end{align}
Certainly, the above classification is closely related to the measures
$\mu$ and $\nu$. For instance, in the DN- and NN-cases,
one requires that $\mu(0,N)<\infty$ and $\mu (-M, N)<\infty$, respectively.
Otherwise, one gets a trivial result
as can be seen from Theorem \ref{t2}{} below.

To state the main result of the paper, we need some assumptions.
In the NN-case (i.e., the ergodic one), we technically assume
that $a$ and $b$ are continuous on $(-M, N)$.
For $\lz^{\text{\rm DN}}$ and $\lz^{\text{\rm NN}}$,
we allow the process to be explosive since the maximal domain is
adopted in definition of $\lz^{\text{\rm DN}}$ and $\lz^{\text{\rm NN}}$.
But for $\lz^{\text{\rm ND}}$ and $\lz^{\text{\rm DD}}$, we are working for the
minimal process (using the minimal domain) only, assuming that $\mu$ and $\nu$
are locally finite.

\thm [Basic estimates\;\ct{cmf10}]\lb{t2}{\cms
Under the assumptions just mentioned, corresponding to each $\#$-case, we have
\be\big(\kz^{\#}\big)^{-1}/4\le \lz^{\#}=\vz_{\max}\le \big(\kz^{\#}\big)^{-1},\de
where
\bg{align}
\big(\kz^{\text{\rm NN}}\big)^{-1}&= \inf_{x<y}
\big[\mu (-M, x)^{-1} + \mu (y, N)^{-1}\big] \nu (x, y)^{-1},\lb{12}\\
\big(\kz^{\text{\rm DD}}\big)^{-1}&= \inf_{x<y}
\big[\nu (-M, x)^{-1} + \nu (y, N)^{-1}\big] \mu (x, y)^{-1},\lb{13}\\
\kz^{\text{\rm DN}}&= \sup_{x\in (0,\, N)} \nu (0, x)\,  \mu (x, N)\lb{14}\\
\kz^{\text{\rm ND}}&= \sup_{x\in (0,\, N)} \mu (0, x)\,  \nu (x, N).\lb{15}
\end{align}}
{\cms In particular, $\lz^{\#}>0$ iff $\kz^{\#}<\infty$.}
\dethm

In each case, the principal eigenvalue is controlled from above and
below by a constant $\kappa^{\#}$ up to a factor 4 which is
universal. Among these cases, the hardest one is the ergodic case.
It may be helpful for the reader to show how to write down
$\kappa^{\text{\rm NN}}$ step by step. \bg{itemize}
\setlength{\itemsep}{-0.8ex}
\item   We need two parameters, say $x$ and $y$ with $x<y$. The state space is then divided by $x$ and $y$
into three parts: the left-hand part $(-M, x)$, the right-hand part $(y, N)$,
and the middle one $(x, y)$.
\item    Measure the left-hand and the right-hand subintervals by $\mu$ and the middle one by $\nu$, respectively:
     $$\kz=\kz^{\text{\rm NN}}:\qqd\qqd \mu (-M, x)\qqd \mu (y, N) \qqd \nu (x, y).$$
\item    Make inverse everywhere:
     $$\kz^{-1}:\qqd\qqd \mu (-M, x)^{-1}\qqd \mu (y, N)^{-1} \qqd \nu (x, y)^{-1}.$$
\item   Finally, summing up the first two terms and making infimum with respect
to $x<y$, we get the answer.
\end{itemize}
Every step is quite natural except the second one: why we use $\mu$
but not $\nu$ in the first two terms? This is because we are in the
ergodic case, $\mu$ is a finite measure. If $\mu$ is replaced by
$\nu$, since $\nu(-\infty, x)$ and $\nu(y, \infty)$ are
infinite when $M$, $N=\infty$, one would get zero for these terms
and so the quantity is trivial. A sensitive point here is that we
use plus, rather than maximum in the last step. Otherwise, even
though the resulting bounds are equivalent to ours but it then would
produce a factor 8 rather than 4 as we expected. We have thus
completed the first, the most important quantity $\kappa^{\text{\rm
NN}}$. To get $\kappa^{\text{\rm DD}}$, simply apply the rule:
exchanging the codes D and N simultaneously in $\kz^{\#}$ leads to
the exchange of the measures $\mu$ and $\nu$ in the formula. Let us
now examine (\ref{13}) more carefully. When $N=\infty$ and $\nu (y,
\infty)=\infty$, the second term in the sum of (\ref{13})
disappeared. In other words, the boundary condition D on the right
endpoint is replaced by N. Then the variable $y$ is free and so can
be removed. Therefore we obtain formula (\ref{14}). We remark
however that the relation between $\lz^{\text{\rm DN}}$ and
$\kappa^{\text{\rm DN}}$ remains the same even if $\nu(y,
\infty)<\infty$.  From (\ref{14}), using again our rule, we obtain
(\ref{15}). We mention that (\ref{15}) can be formally obtained from
(\ref{12}) by removing the second term in the sum. Actually,
(\ref{15}) is formally a reverse of (\ref{14}), and so is somehow an
easy consequence of (\ref{14}).

\subsection{Short review on the known results}

It is the position to say a little about the history of the topic.
Clearly, we are in the typical situation of the Sturm--Liouville
eigenvalue problem (1836-1837). From which, we learn the general
properties of the eigenfunction: the existence and uniqueness, the
zeros of the eigenfunction, and so on. Except some very specific
cases, the problem is usually not solvable analytically. This leads
to the theory of special functions used widely in sciences. The
estimation of the principal eigenvalues is usually not included in
the Sturm--Liouville theory but is studied in harmonic analysis
(especially for $\lz^{\text{\rm DN}}$). To see this, rewrite
(\ref{11}) as the {\it Poincar\'e inequality}
$$\lz^{\text{\rm DN}} \|f\|^2\le D(f),\qqd f(0)=0.$$
More general, we have {\it Hardy's inequality}
$$\|f\|_{L^p(\mu)}^p \le A_p\int_{-M}^N {|f'|}^p e^C,\qqd f(0)=0,\; p> 1$$
where $A_p$ denotes the optimal constant in the inequality. Certainly,
$A_2= \big(\lz^{\text{\rm DN}}\big)^{-1}$. This was initialed,
for the specific operator $L=x^2 \d^2/\d x^2$, by G.H. Hardy \ct{hgh20} in 1920,
motivated from a theorem of Hilbert on double series. To which,
several famous mathematicians (H. Weyl, F.W. Wiener, I. Schur, et al.)
were involved.
After a half-century, the basic estimates in the DN-case were
finally obtained by several mathematicians, for instance B. Muckenhoupt (1972).
The reason should be now clear why (\ref{14}) can be
so famous in the history. The estimate of $\lz^{\text{\rm ND}}$
was given in Maz'ya (1985). In the DD-case, the problem was begun by P. Gurka
\ct{gp89} around 1989 and then improved in
the book by Opic and Kufner (1990) with a factor $\approx 22$. In terms of a splitting technique, the NN-case
can be reduced to the Muckenhoupt's estimate with a factor 8, as shown by Miclo (1999) in the context
of birth--death processes. A better estimate can be
done in terms of variational formulas given in \rf{cmf00}{Theorem 3.3}. It is surprising that in the more complicated DD-
and NN-cases, by adding one more parameter only, we can still obtain a compact expression
(\ref{12}) and (\ref{13}). Note that these two formulas have the following advantage:
 the left- and the right-hand
 parts are symmetric; the cases having finite or infinite intervals are unified together
  without using the splitting technique.

\subsection{Motivation and application}

Here is a quick overview of our motivation and application of the
study on this topic. Consider the $\fz^4$-model on the
$d$-dimensional lattice ${\mathbb Z}^d$. At each site $i$, there is
a one-dimensional diffusion with operator $L_i=\d^2/\d x_i^2
-u'(x_i)\d/\d x_i$, where $u_i(x_i)=x_i^4-\bz x_i^2$ having a
parameter $\bz\ge 0$. Between the nearest neighbors $i$ and $j$ in
${\mathbb Z}^d$, there is an interaction. That is, we have an
interaction potential
 $H(x)=-J\sum_{\langle i j\rangle} x_i x_j $ with parameter $J\ge 0$.
For each finite box $\llz$ (denoted by $\llz\Subset {\mathbb Z}^d$) and $\oz\in {\mathbb R}^{{\mathbb Z}^d}$,
let $H_{\llz}^{\oz}$ denote the conditional Hamiltonian (which acts on those $x$: $x_k=\oz_k$ for all
$k\notin \llz$). Then, we have a local operator
$$L_{\llz}^{\oz}=\sum_{i\in \llz}\big[\partial_{ii}-\partial_i(u+H_{\Lambda}^{\omega})\partial_i\big].$$
 We proved that the first non-trivial
eigenvalue $\lambda_1^{\bz} \big(\Lambda, \omega \big)$ (as well as the logarithmic Sobolev constant
$\sz^{\bz} \big(\Lambda, \omega \big)$ which is not touched here) of $L_{\llz}^{\oz}$ is approximately
$\exp[-\bz^2/4]- 4 d J$ uniformly with respect to the boxes $\llz$
and the boundaries $\oz$. The leading
rate $\bz^2/4$ is exact which is the only one we have ever known up to now for a continuous model.

\thm[\ct{cmf08}] {\cms For the $\fz^4$-model given above, we have
$$\begin{aligned}
&\inf_{\Lambda\Subset {\mathbb Z}^d}\inf_{\omega\in {\mathbb
R}^{\mathbb Z^d}} \lambda_1^{\bz} \big(\Lambda, \omega \big)
\!\approx\! \inf_{\Lambda\Subset
{\mathbb Z}^d}\inf_{\omega\in {\mathbb R}^{\mathbb Z^d}} \sz^{\bz}
\big(\Lambda, \omega \big) \!\approx\! \exp\big[\!-\!\bz^2/4\!-\!c\log
\bz\big]-4dJ,\end{aligned}$$
where $c\in [1, 2]$. See Figure 1.}
\begin{center}{\includegraphics{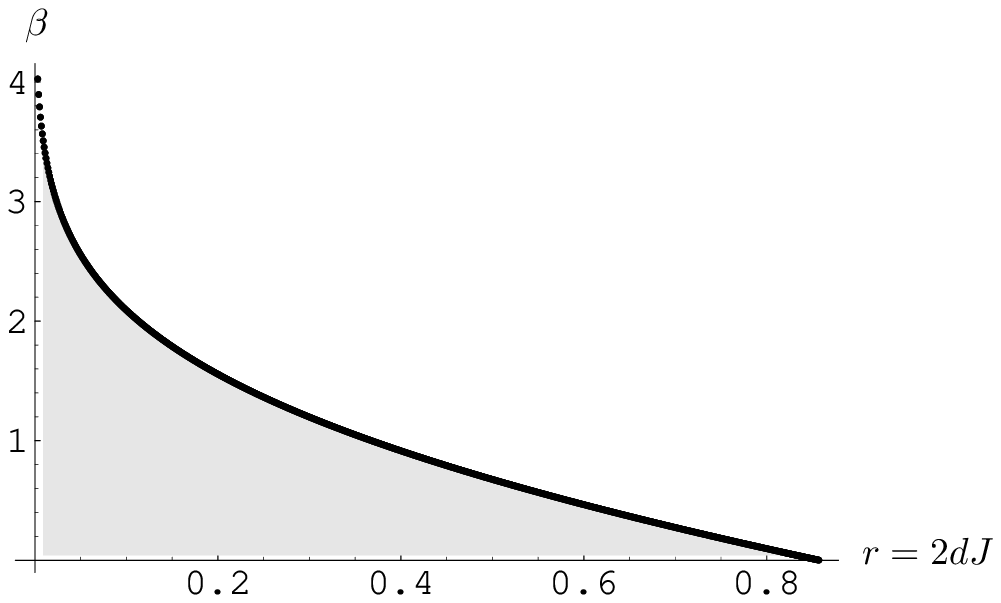}
\vspace{-6truecm}
$$\aligned
&\lz_1^{\bz}, \; \sz^{\bz}:\\
&\exp\big[\!-\!\bz^2/4\!-\!c\log
\bz\big]-2r\\
&\qqd\qd c=c(\bz)\in [1, 2]\endaligned$$

\vspace{3.5 truecm}

{\bf Figure 1}\qd\rm Phase transition of the $\fz^4$ model}\end{center}
\dethm

The figure says that in the gray region, the system has a positive
principal eigenvalue and so is ergodic; but in the region which is a
little away above the curve, the eigenvalue vanishes. The picture
exhibits a phase transition. The key to prove Theorem \ref{t2} is a
deep understanding about the one-dimensional case. Having
one-dimensional result at hand, as far as we know, there are at
least three different ways to go to the higher or even infinite
dimensions:  the conditional technique used in \ct{cmf08}; the
coupling method explained in \rf{cmf05a}{Chapter 2}; and some
suitable comparison which is often used in studying the stability
rate of interacting particle systems. This explains our original
motivation and shows the value of a sharp estimate for the leading
eigenvalue in dimension one. The application of the present result
to this model should be clear now.

\section{Sketch of the proof}\lb{s2}

The hardest part of Theorem \ref{t2} is the assertion for
$\lz^{\text{NN}}$. Here we sketch its proof. Meanwhile, the proof
for $\lz^{\text{DD}}$ is also sketched. The proof for the first
assertion consists mainly of three steps by using three methods: the
coupling method, the dual method, and the capacitary method.

\subsection{Coupling method}

The next result was proved by using the coupling technique.

\thm[Chen and Wang (1997)] {\cms For the operator $L$ on $(0,
\infty)$ with reflection at $0$, we have}
\begin{align}
\lz_1&=\lz^{\text{\rm NN}}\ge \sup_{f\in {\scr F}}\,
\inf_{x>0}\bigg[-b' - \frac{a f'' + (a'+b)f'}{f}\bigg](x), \\
{\scr F}&=\big\{f\in {\scr C}^2(0, \infty): f(0)=0,\; f|_{(0,\, \infty)}>0\big\}.
\end{align}
{\cms Actually, the equality sign holds once the eigenfunction
ot $\lz_1$ belongs to $C^3$}.
\dethm

We now rewrite the above formula in terms of an operator, Schr\"odinger operator $L_S$.
\bg{align}
{\lz_1}&= \sup_{f\in {\scr F}}\,\inf_{x>0}\bigg[-b' - \frac{a f'' + (a'+b)f'}{f}\bigg](x) \\
&= \sup_{f\in {\scr F}}\,\inf_{x>0} \bigg(-\frac{L_S\, f}{f}\bigg)(x)=:{\lz_S},\lb{19}\\
L_S&=a(x)\frac{\d^2}{\d x^2}+ \big(a'(x)+b(x)\big)\frac{\d}{\d x}+ {b'(x)}.
\end{align}

The original condition $\pi (f)=0$ in the definition of
$\lz^{\text{NN}}$ means that $f$ has to change its sign. Note that
$f$ is regarded as a mimic of the eigenfunction $g$. The difficulty
is that we do not know where the zero-point of $g$ is located. In
the new formula (\ref{19}), the zero-point of $f\in {\scr F}$ is
fixed at the boundary $0$, the function is positive inside of the
interval. This is the advantage of formula (\ref{19}). Now, a new
problem appears: there is an additional potential term $b'(x)$.
Since $b'(x)$ can be positive, the operator $L_S$ is Schr\"odinger
but may not be an elliptic operator with killing. Up to now, we are
still unable to handle with general Schr\"odinger operator (even
with killing one), but at the moment, the potential term is very
specific so it gives a hope to go further.

\subsection{Dual method}

To overcome the difficulty just mentioned, the idea is a use of duality. The dual now we adopted is very simple:
just an exchange of the two measures $\mu$ and $\nu$. Recall that the original
operator is $L=\frac{\d}{\d \mu} \,\frac{\d}{\d \nu}$ by (\ref{7}). Hence the dual operator takes the following form
\bg{align}
L^*&=\frac{\d}{\d \mu^*} \,\frac{\d}{\d \nu^*}={\frac{\d}{\d \nu} \,\frac{\d}{\d \mu}},\\
L^*&= a(x)\frac{\d^2}{\d x^2}+ \big(a'(x)-b(x)\big)\frac{\d}{\d x},\qqd
x\in (0, \infty).
\end{align}
This dual goes back to Siegmund (1976) and Cox \& R\"osler (1983) (in which the
probabilistic meaning of this duality was explained), as an analog
of the duality for birth--death process (cf. \ct{cmf10} for more
details and original references). It is now a simple matter to check
that the dual operator is a similar transform of the Schr\"odinger
one \be L^*=e^C L_S e^{-C}.\de
This implies that
$$-\frac{L_S\, f}{f} =-\frac{L^* f^*}{f^*},$$
where $f^*:=e^C f$ is one-to-one from ${\scr F}$ into itself.
Therefore, we have
$$\lz_S=\sup_{f\in {\scr F}}\,\inf_{x>0} \frac{-L_S\,
f}{f}(x)=\sup_{f^*\in {\scr F}}\,\inf_{x>0} \frac{-L^* f^*}{f^*}(x)=\lz^{*\text{\rm DD}},$$
where the last equality is the so-called {\it Barta's equality}.

we have thus obtained the following identity.

\prp \lb{t22} $\lz_1=\lz_S=\lz^{*\text{\rm DD}}$.
\deprp

Actually, we have a more general conclusion that $L_S$ and $L^*$ are isospectral
from $L^2\big(e^C\d x\big)$ to
$L^2\big(e^{-C}\d x\big)$. This is because of
$$\int e^C f\, L_S\, g=\int e^{-C}(e^C f) \big(e^C L_S e^{-C}\big)(e^C g)
=\int e^{-C} f^* L^* g^*,$$
and $L_S$ and $L^*$ have a common core. But $L$ on $L^2(\mu)$ and its dual $L^*$ on $L^2\big(e^{-C}\d x\big)$
are clearly not isospectral.

The rule mentioned in the remark after Theorem \ref{t2},
and used to deduced (\ref{13}) from  (\ref{12}),
comes from this duality.
Nevertheless, it remains to compute $\lz^{\text{DD}}$ for the dual operator.

\subsection{Capacitary method}

To compute $\lz^{\text{DD}}$, we need a general result which comes
from a different direction to generalize the Hardy-type inequalities.
In contract to what we have talked so far, this time we extend
the inequalities to the higher dimensional situation.
This leads to a use of the capacity since in the higher
dimensions, the boun\-dary may be very complicated. After a great effort by many
mathematicians (see for instance Maz'ya 1985; Hasson 1979; Vondra\v cek 1996;
Fukushima \& Uemura 2003; and \ct{cmf05b}), we have finally the following result.

\thm {\cms For a regular transient Dirichlet form  $(D, {\scr D}(D))$
with locally compact state space $(E, {\scr E})$,
the optimal constant $A_{\mathbb B}$ in the Poincar{\cms\'e}-type
inequality
$$\big\|f^2\big\|_{\mathbb B}\le A_{\mathbb B}\, D(f),\qqd f\in {\scr C}_K^\infty (E)$$
satisfies
$ B_{\mathbb B}\le A_{\mathbb B}\le 4 B_{\mathbb B}, $
where $\|\cdot\|_{\mathbb B}$ is the norm in a normed linear space ${\mathbb B}$ and}
$$B_{\mathbb B}=\sup_{\text{\rm compact}\, K} {\text{\rm Cap}(K)}^{-1}{{\|\mathbbold{1}_K\|_{\mathbb B}}}. $$
\dethm

The space ${\mathbb B}$ can be very general, for instance $L^p(\mu)\,(p\ge 1)$ or the
Orlicz spaces. In the present context,
$D(f)= \int_{-M}^N {f'}^2 e^C$,
${\scr D}(D)$ is the closure of ${\scr C}_K^{\infty}(-M, N)$
with respect to the norm $\|\cdot\|_D$: $\|f\|_D^2\!=\!\|f\|^2\!+\!D(f)$, and
$$\text{\rm Cap}(K)=\inf\big\{D(f): f\in {\scr C}_K^{\infty}(-M, N),\, f|_K\ge 1\big\}.$$

Note that we have the
universal factor 4 here and the isoperimetric constant
$B_{\mathbb B}$ has a very compact form.
We now need to compute the capacity only. The problem is that the capacity
is usually not computable explicitly. For instance, at the moment, I do not
know how to compute it for Schr\"odinger operators even for the
elliptic operators having killings. Very lucky,
we are able to compute the capacity for the one-dimensional
elliptic operators. The result has a simple expression:
$$B_{{\mathbb B}}=\sup_{-M< x < y <N}
\big[\nu(-M, x)^{-1} + \nu(y, N)^{-1}\big]^{-1}\|\mathbbold{1}_{(x,\, y)}\|_{{\mathbb B}}.$$
It looks strange to have double inverse here. So, making inverse in both sides, we get
$$B_{{\mathbb B}}^{-1}=\inf_{-M< x < y <N}
\big[\nu(-M, x)^{-1} + \nu(y, N)^{-1}\big]\,\|\mathbbold{1}_{(x, y)}\|_{{\mathbb
B}}^{-1}.$$
Applying this result to ${\mathbb B}=L^1(\mu)$, we obtain
the solution to the DD-case: $\lz^{\text{\rm DD}}=A_{L^1(\mu)}^{-1}$ and
$$\big(\kz^{\text{\rm DD}}\big)^{-1}
=B_{L^1(\mu)}^{-1}=\inf_{-M< x < y <N}\big[\nu(-M, x)^{-1} + \nu(y, N)^{-1}\big]\,
{\mu (x, y)^{-1}}.$$

\subsection{The final step}

Applying the last result to the dual process and using Proposition \ref{t22}, we have
not only
$$\big(\kz^{*\text{\rm DD}}\big)^{-1}/4
\le \lz^{\text{\rm NN}}=\lz_S= \lz^{*\text{\rm DD}}
\le \big(\kz^{*\text{\rm DD}}\big)^{-1},$$
but also
$$\aligned
{\big(\kz^{*\text{\rm DD}}\big)^{-1}}\!\!\!&\!=\!
\inf_{x < y}\big[\nu^*\!(-M, x)^{-1}\!\! +\! \nu^*\!(y, N)^{-1}\big]\,
\mu^*\! (x, y)^{-1}\\
&=\!
\inf_{x < y}\big[\mu(-M, x)^{-1}\!\! +\! \mu(y, N)^{-1}\big]\,
\nu(x, y)^{-1}\\
&={\big(\kz^{\text{\rm NN}}\big)^{-1}}.
\endaligned$$
This finishes the proof of the main assertion of Theorem \ref{t2}.

\subsection{Summery of the proof}

Here is the summery of our proof. First, by a change of the topology,
we reduce the study on $\vz_{\max}$ to $\lz^{\text{NN}}$. Then, by coupling,
we reduce $\lz^{\text{NN}}$ to $\lz_S$. Next, by duality, we reduce $\lz_S $
to $\lz^{*\text{DD}}$. We use capacitary method to compute $\lz^{*\text{DD}}$.
Finally, we use duality again to come back to $\lz^{\text{NN}}$. Recall that our
original purpose is using $\lz_1=\lz^{\text{NN}}$ to study the phase transition,
a basic topic in the study on interacting particle systems (abbrev. IPS).
It is very interesting that we now have an opposite interaction. We use the main tools
(coupling and duality) developed in the study on IPS to investigate a very classical problem
and produce an interesting result.

\section{Improvements}\lb{s3}

The basic estimates given in Theorem \ref{t2}{} can be further improved.
For half-line at least, we have actually
 an approximating procedure for each of the principal
eigenvalues. Refer to \ct{cmf05a, cmf10} and references therein.
Moreover, one may approach the whole line by half-lines.
Here we consider an additional method but concentrate
on $\lz^{\text{DD}}$ and $\lz^{\text{NN}}$ only.
As will be seen soon, the resulting bounds are much more complicated,
less simple and less symmetry, than those given in Theorem \ref{t2}.

Let us begin with a simper but effective result.

\prp\lb{t031} {\cms We have
$$\lz^{\text{\rm DD}}
\le \big({\bar\kz}^{\text{\rm DD}}\big)^{-1}\le
\big(\kz^{\text{\rm DD}}\big)^{-1}$$
and
$$\lz^{\text{\rm NN}}
\le \big({\bar\kz}^{\text{\rm NN}}\big)^{-1}
\le \big(\kz^{\text{\rm NN}}\big)^{-1}, $$
where}
\begin{align*}
&\big({\bar\kz}^{\text{\rm DD}}\big)^{-1}= \inf_{x<y}\Big({\nu(-M, x)}^{-1}+{\nu[y, N)}^{-1}\Big)\times\\
&\qd \times\bigg\{ \mu(x, y)+\!\!\int_{-M}^x\mu (\d z)\bigg[1-\frac{\nu(z,x)}{\nu (-M, x)}\bigg]^2
\! +\!\!\int_y^N \!\!\mu (\d z)\bigg[1- \frac{\nu(y, z)}{\nu(y, N)}\bigg]^2\bigg\}^{-1}\!,\\
&\big({\bar\kz}^{\text{\rm NN}}\big)^{-1}= \inf_{x<y}\Big({\mu(-M, x)}^{-1}+{\mu[y, N)}^{-1}\Big)\times\\
&\qd \times\bigg\{ \nu(x, y)+\!\!\int_{-M}^x\nu (\d z)\bigg[1-\frac{\mu(z,x)}{\mu (-M, x)}\bigg]^2
\! +\!\!\int_y^N \!\!\nu (\d z)\bigg[1- \frac{\mu(y, z)}{\mu(y, N)}\bigg]^2\bigg\}^{-1}\!.
\end{align*}
\deprp

Note that if $\nu(-M, N)<\infty$ which is not assumed in
Proposition \ref{t031}, then the last two terms
in $\{\cdots\}$ in the expression
of $\big({\bar\kz}^{\text{\rm DD}}\big)^{-1}$ can be written as
$$\nu(-M, x)^{-2}\!\int_{-M}^x\mu (\d z) \nu(-M, z)^2
\! +\nu(y, N)^{-2}\!\int_y^N \!\!\mu (\d z) \nu(z, N)^2.$$
Otherwise, this expression may be meaningless.
Similar comment is meaningful for
$\big({\bar\kz}^{\text{\rm NN}}\big)^{-1}$.

\prf  Fix $x<y$. Applying
$\lz^{\text{\rm DD}}\le {D(f)}/{\mu\big(f^2\big)}$
to the test function
$$f(z)=
\begin{cases}
{\displaystyle\frac{\nu(y, N)}{\nu (-M, x)}}\, \nu (-M, z\wedge x),\qd & z\le y\\
\nu (z, N), & z\ge y.
\end{cases}
$$
we obtain $\lz^{\text{\rm DD}}
\le \big({\bar\kz}^{\text{\rm DD}}\big)^{-1}$.
By duality, we obtain the assertion for ${\bar\kz}^{\text{\rm NN}}$.
Refer to the remark after the proof of \rf{cmf10}{Theorem 8.2}
for more details.
\deprf

To improve the lower estimate Theorem \ref{t2}{}, we need more work.
For a given $f\in {\scr C}(-M, N)$ with $f|_{(-M, N)}>0$, define
\begin{align}
h^-(z)&\!=\!h_f^-(z)\!=\!\nu\Big(\mu \big(\mathbbold{1}_{(\cdot,\, \uz)}f\big) \mathbbold{1}_{(-M,\, z)}\Big)
\!=\!\!\int_{-M}^{z}\! e^{-C(x)}\d x \int_x^{\uz}\!\frac{e^C f}{a},\;
 z\le \uz,\lb{3-25} \\
h^+(z)&\!=\!h_f^+(z)\!=\!\nu\Big(\mu \big(\mathbbold{1}_{(\uz,\, \cdot)}f\big) \mathbbold{1}_{(z,\, N)}\Big)
\!=\!\!\int_{z}^{N}\! e^{-C(x)}\d x \int_{\uz}^x\!\frac{e^C f}{a},\qd z>\uz, \lb{3-26}
\end{align}
i.e. (by exchanging the order of the integrals),
$$\aligned
h^-(z)&\!=\!\mu\big(f \nu(-M, \cdot \wedge z)\big)\!\!=\!\!
\mu\Big(f\nu (-M, \cdot) \mathbbold{1}_{(-M, z)}\Big)\!\!+\!\mu\Big(\!f\mathbbold{1}_{(z, \uz)}\!\Big) \nu(-M, z),\, z\le \uz, \\
h^+(z)&=\mu\big(f \nu(\cdot \vee z, N)\big)=\!
\mu\Big(f\nu (\cdot, N) \mathbbold{1}_{(z,\, N)}\Big)\!+\mu\Big(f\mathbbold{1}_{(\uz,\, z)}\Big) \nu(z, N),\qd\;\;\, z>\uz,
\endaligned$$
where $x \wedge y=\min\{x, y\}$, $x \vee y=\max\{x, y\}$, and
$\uz=\uz (f)\in (-M, N)$ is the unique root of the equation:
$$h^-(\uz)= h^+(\uz)$$
provided $h_f^{\pm}<\infty$. Next, define $I\!I^{\pm}(f)= h^{\pm}/f$.

\thm[Variational formula]\lb{t032}{\cms Let $a$ and $b$ be continuous and $a>0$ on $(-M, N)$.
\begin{itemize}
\item [(1)] Assume that $\nu(-M, N)<\infty$.
Using the notation above, we have
\be\lz^{\text{\rm DD}}= \sup_{f\in {\scr C}_+}\Big[\inf_{z\in (-M, \uz)}I\!I^-(f)(z)^{-1}\Big]
\bigwedge \Big[\inf_{z\in (\uz, N)}I\!I^+(f)(z)^{-1}\Big],\lb{027}\de
where ${\scr C}_+=\{f\in {\scr C}(-M, N): f>0 \;\text{\cms on}\; (-M, N) \}$.
\item [(2)] Assume that $\mu(-M, N)<\infty$. Then (\ref{027})
holds replacing $\lz^{\text{\rm DD}}$ by $\lz^{\text{\rm NN}}$ provided in definition of $h^{\pm}$, $\mu$ and $\nu$ are exchanged.
\end{itemize}
}
\dethm

\prf By duality, it suffices to prove the first assertion.

(a) Without loss of generality, assume that $h_f^{\pm}<\infty$.
Otherwise, the assertion is trivial. First, we prove ``$\ge$''. Let
$$h(z)=
\begin{cases}
h^-(z),\qqd z\le \uz,  \\
h^+(z),\qqd z>\uz,
\end{cases}$$
Clearly, $h|_{(-M, N)}>0$ and
$h\in {\scr C}(-M, N)$ in view of definition of $\uz$. Next, note that
$$\aligned
&{h^-}'(x)= e^{-C(x)}\int_x^{\uz} \frac{e^C}{a}f,\qd
{h^-}''(x)= e^{-C(x)}\bigg[-\frac{b}{a}
\int_x^{\uz} \frac{e^C}{a}f- \frac{e^C}{a}f\bigg],\qd x<\uz;\\
&{h^+}'(x)= -e^{-C(x)}\int_{\uz}^x \frac{e^C}{a}f,\qd
{h^+}''(x)= e^{-C(x)}\bigg[\frac{b}{a}
\int_{\uz}^x \frac{e^C}{a}f- \frac{e^C}{a}f\bigg],\qd x>\uz.
\endaligned$$
Obviously, $h'(\uz\pm 0)=0$. Since $a$, $b$ and $f$ are continuous and $a>0$ on $(-M, N)$,
we also have $h''(\uz+0)=h''(\uz-0)$ and so
$h\in {\scr C}^2(-M, N)$. Therefore, by Barta's equality, we have
$$\aligned
\lz^{\text{\rm DD}}&=\sup_{g\in {\scr F}}\,\inf_{z\in (-M, N)}\frac{-L g}{g}(z)\\
&\ge \inf_{z\in (-M, N)}\frac{-L h}{h}(z)\\
&=\Big[\inf_{z\in (-M, \uz)}\frac{-L h^-}{h^-}(z)\Big]
\bigwedge \Big[\inf_{z\in (\uz, N)}\frac{-L h^+}{h^+}(z)\Big].
\endaligned$$
Now, by (\ref{7}), required assertion follows by a simple computation.

(b) Next, we show that the equality sign in (\ref{027}) holds.
The assertion becomes trivial if $\lz^{\text{\rm DD}}=0$. Otherwise,
the eigenfunction $g$ of $\lz^{\text{\rm DD}}$ should be unimodal (which seems known in the Sturm--Liouville theory and is proved in the
discrete context \rf{cmf10}{Proposition 7.14}. Actually, the
discrete case is even more complex since the eigenfunction can
be a simple echelon, not necessarily unimodal).
By setting $f=g$ and $\uz$ to be the maximum point of $g$ ($g'(\uz)=0$), it follows that $I\!I^{\pm}(f)^{-1}\equiv \lz^{\text{\rm DD}}$ and hence the equality sign holds.
\deprf

We now introduce a typical application of Theorem \ref{t032}{}. Fix $x<y$. Define
$$f^{x, y}(s)=
\begin{cases} {\displaystyle\sqrt{\frac{\nu(y, N)}{\nu (-M, x)}\, \nu (-M, s\wedge x)}},\qd & s\le y\\
\sqrt{\nu (s, N)}, & s\ge y
\end{cases}
$$
and set
$${\underline\kz}^{\text{\rm DD}}=\inf_{x< y}
\Big[\sup_{z\in (-M, \uz)}I\!I^-(f^{x, y})(z)\Big]
\bigvee \Big[\sup_{z\in (\uz, N)}I\!I^+(f^{x, y})(z)\Big].
 $$
By exchanging $\mu$ and $\nu$, we obtain ${\underline\kz}^{\text{\rm NN}}$.
Now, by Theorem \ref{t032}{},
we have the following result.

\crl\lb{t033} {\cms Under assumptions of Theorem \ref{t032}{}, we have
$$\lz^{\text{\rm DD}}\ge \big({\underline\kz}^{\text{\rm DD}}\big)^{-1}
\qd\text{\cms and}\qd \lz^{\text{\rm NN}}\ge \big({\underline\kz}^{\text{\rm NN}}\big)^{-1}.$$
}
\decrl

We remark that the assumption in part (1) of Theorem \ref{t032} is
necessary for DD-case (cf. (\ref{12})). Recall that (\ref{027}) is
a complete variational formula for the lower estimates of $\lz^{\text{\rm DD}}$. Starting
at $f_1=f$ used in Corollary \ref{t033}, replacing $f$ and $h$ used in Theorem \ref{t032} by $f_{n-1}$ and $f_n$, respectively, we obtain an approximating procedure from below for $\lz^{\text{\rm DD}}$.
Dually, we can obtain a variational formula for the upper estimates
of $\lz^{\text{\rm DD}}$ and an approximating procedure from above.
Here we omit all of the details.
The same remark is meaningful for $\lz^{\text{\rm NN}}$, which is especially interesting
since here we do not use the property that $\mu (f)=0$ for the test function
$f$. The new difficulty of (\ref{027}) is that $\uz(f)$ may not be
computable analytically. This costs a question to prove that
${\underline\kz}^{\text{\rm DD}}\le 4 {\kz}^{\text{\rm DD}}$
which should be true in view of our knowledge on the half-line,
and is illustrated by examples in the next section.
It is noticeable that the method works for the whole line
and the use of $\uz(f)$ is essentially
different from what used in the splitting technique.
Finally, we mention that the method used here is meaningful for
birth--death processes, refer to \rf{cmf10}{Lemma 7.12}.

For convenience in practice, we express $h^{\pm}$ used in Corollary
\ref{t033} more explicitly. Let $\nu_-(s)=\nu (-M, s)$ and
$\nu_+(s)=\nu(s, N)$ for simplicity. Then
\be f(s)=f^{x, y}(s)=
\begin{cases}
\sqrt{\nu_+(y)\nu_-(s)}\big/\sqrt{ \nu_-(x)}, &\qd s\le x\\
\sqrt{\nu_+ (y)}, &\qd x\le s\le y\\
\sqrt{\nu_+(s)}, &\qd s\ge y,
\end{cases}\lb{3-28}\de
and
\begin{align}
h^-(z)&= \mu\Big(f\nu_- \mathbbold{1}_{(-M,\, z)}\Big)
+\nu_- (z)\, \mu \big(f \mathbbold{1}_{(z,\, \uz)}\big),\qqd z\le \uz,\lb{3-29}\\
h^+(z)&= \mu\Big(f\nu_+ \mathbbold{1}_{(z,\, N)}\Big)
+\nu_+ (z)\, \mu \big(f \mathbbold{1}_{(\uz,\, z)}\big),\qqd\qd z\ge \uz.\lb{3-30}
\end{align}
We now consider the typical case that $ \uz \in [x, y]$.
Then,
$$\aligned
h^-(\uz)&=\sqrt{\frac{\nu_+(y)}{\nu_-(x)}}\,
\mu\Big(\nu_-^{3/2} \mathbbold{1}_{(-M,\, x)}\Big)
+\sqrt{\nu_+(y)}\,
\mu\Big(\nu_- \mathbbold{1}_{(x,\, \uz)}\Big),\\
h^+(\uz)&=
\mu\Big(\nu_+^{3/2} \mathbbold{1}_{(y,\, N)}\Big)
+\sqrt{\nu_+(y)}\,
\mu\Big(\nu_+ \mathbbold{1}_{(\uz,\, y)}\Big).
\endaligned$$
Hence the equation $h^-(\uz)=h^+(\uz)$ becomes
\begin{align}
&\frac{1}{\sqrt{\nu_-(x)}}\mu\Big(\nu_-^{3/2} \mathbbold{1}_{(-M,\, x)}\Big)
+\mu\Big(\nu_- \mathbbold{1}_{(x, \uz)}\Big)\nonumber\\
&\qd
= \frac{1}{\sqrt{\nu_+(y)}}
\mu\Big(\nu_+^{3/2} \mathbbold{1}_{(y,\, N)}\Big)
+\mu\Big(\nu_+ \mathbbold{1}_{(\uz,\, y)}\Big),\qqd \uz\in [x, y].\lb{3-31}\end{align}
Furthermore, by some computations, we obtain the ratio $h^{\pm}/f^{x,y}$
as follows. We have for $z$: $z\le x\le \uz \le y$ that
\begin{align}
I\!I^-\big(f^{x, y}\big)(z)=&\frac{1}{\sqrt{\nu_-(z)}}\,
\mu\Big(\nu_-^{3/2} \mathbbold{1}_{(-M,\, z)}\Big)
+ \sqrt{{\nu_-(z)}}\,
\mu\Big(\sqrt{\nu_-}\, \mathbbold{1}_{(z,\, x)}\Big)\nonumber\\
&+\sqrt{\nu_-(z)\nu_-(x)}\,\mu(x, \uz),\lb{3-32} \end{align}
and for $z$: $z\ge y\ge \uz$ that
\begin{align} I\!I^+\big(f^{x, y}\big)(z)
=&\frac{1}{\sqrt{\nu_+(z)}}\mu\Big(\nu_+^{3/2} \mathbbold{1}_{(z,\, N)}\Big)
+ \sqrt{\nu_+(z)}\,\mu\Big(\sqrt{\nu_+}\, \mathbbold{1}_{(y,\, z)}\Big)
\nonumber\\
&+ \sqrt{\nu_+(z)\nu_+(y)}\,\mu(\uz, y).\lb{3-33}\end{align}
Note that by (\ref{3-25}) and (\ref{3-26}), $h^-$ is increasing on $[x, \uz]$ and
$h^+$ is decreasing on $[\uz, y]$. Since $f^{x, y}$ is a constant on
$[x, y]$, it follows that
$$\max_{z\in [x, \uz]} \frac{h^-(z)}{f^{x, y}(z)}= \frac{h^-(\uz)}{f^{x, y}(x)}
\qd \text{and}\qd
\max_{z\in [\uz, y]} \frac{h^+(z)}{f^{x, y}(z)}= \frac{h^+(\uz)}{f^{x, y}(x)}.$$
By assumption, $h^-(\uz)=h^+(\uz)$. Hence
\begin{align}
\max_{z\in [x, \uz]} I\!I^-\big(f^{x, y}\big)(z)
&= \max_{z\in [\uz, y]} I\!I^+\big(f^{x, y}\big)(z)
=\frac{h^-(\uz)}{f^{x, y}(x)}\nonumber\\
&=\frac{1}{\sqrt{\nu_-(x)}}\,
\mu\Big(\nu_-^{3/2} \mathbbold{1}_{(-M,\, x)}\Big)
+ \mu\Big(\nu_-\, \mathbbold{1}_{(x,\, \uz)}\Big).\lb{3-34}
\end{align}
Thus, for computing ${\underline\kz}^{\text{\rm DD}}$, by (\ref{3-32})--(\ref{3-34}), we arrive at
\begin{align}
&\Big[\sup_{z\in (-M,\, \uz)}I\!I^-(f^{x, y})(z)\Big]
\bigvee \Big[\sup_{z\in (\uz, N)}I\!I^+(f^{x, y})(z)\Big]\nonumber\\
&\qd =\sup_{z\in (-M,\, x)}\bigg[\frac{1}{\sqrt{\nu_-(z)}}\,
\mu\Big(\nu_-^{3/2} \mathbbold{1}_{(-M,\, z)}\Big)
+ \sqrt{{\nu_-(z)}}\,
\mu\Big(\sqrt{\nu_-}\, \mathbbold{1}_{(z,\, x)}\Big)\nonumber\\
&\qqd\qqd\qqd\; +\sqrt{\nu_-(z)\nu_-(x)}\,\mu(x,\, \uz)\bigg]\nonumber\\
&\qqd\bigvee \bigg[\frac{1}{\sqrt{\nu_-(x)}}\,
\mu\Big(\nu_-^{3/2} \mathbbold{1}_{(-M,\, x)}\Big)
+ \mu\Big(\nu_-\, \mathbbold{1}_{(x,\, \uz)}\Big)\bigg]\nonumber\\
&\qqd\bigvee\sup_{z\in (y,\, N)} \bigg[\frac{1}{\sqrt{\nu_+(z)}}\mu\Big(\nu_+^{3/2} \mathbbold{1}_{(z,\, N)}\Big)
+ \sqrt{\nu_+(z)}\,\mu\Big(\sqrt{\nu_+}\, \mathbbold{1}_{(y,\, z)}\Big)
\nonumber\\
&\qqd\qqd\qqd\;\;\, + \sqrt{\nu_+(z)\nu_+(y)}\,\mu(\uz, y)\bigg].
\lb{3-35}\end{align}
Finally, let $(x^*, y^*, \uz^*)$ solve equation (\ref{3-31}) and two more equations
modified from (\ref{3-35})
ignoring its left-hand side and replacing the last two
``$\vee$'' with ``$=$''. Then we have
\be
{\underline\kz}^{\text{\rm DD}} =\frac{1}{\sqrt{\nu_-(x^*)}}\,
\mu\Big(\nu_-^{3/2} \mathbbold{1}_{(-M,\, x^*)}\Big)
+ \mu\Big(\nu_-\, \mathbbold{1}_{(x^*,\, \uz^*)}\Big).\lb{3-36}\de

\section{Examples}\lb{s4}

This section illustrates the application of the basic estimates
given in Theorem \ref{t2} and the improvements given in Proposition
\ref{t031} and Corollary \ref{t033}.

\xmp [OU-processes]\lb{t041} {\rm The state space is $\mathbb R$ and the operator
is
$$L=\frac{1}{2}\bigg(\frac{\d^2}{\d x^2}- 2 x \frac{\d}{\d x}\bigg).$$
This is a typical example of the use of special functions.
It has discrete eigenvalues $\lz_n=n$ with eigenfunctions (Hermite polynomials)
$$g_n(x)=(-1)^n e^{x^2} \frac{\d^n}{\d x^n}\big(e^{-x^2}\big),\qqd n\ge 0.$$
Then, we have $\big(\kz^{\text{\rm DD}}\big)^{-1}=\lz_0=0$,
$\lz^{\text{\rm NN}}=\lz_1=1$ with eigenfunction $g(x)=x$. To compute
$\kz^{\text{\rm NN}}$, noting that the operator, the eigenfunction are all
symmetric with respect to $0$ and so does $\kz^{\text{\rm NN}}$,
one can split the whole line into two parts
$(-\infty, 0)$ and $(0, \infty)$ with common Dirichlet boundary at $0$.
This simplifies the computation and leads to
$\big(\kz^{\text{\rm NN}}\big)^{-1}=\big(\kz^{\text{\rm DN}}\big)^{-1}\approx 2.1$.
Note that $g'(x)\equiv 1$ but $\lim_{|x|\to\infty}\big(e^C g'\big)(x)=0$.

For the half-space $(0, \infty)$, as we have just mentioned,
$\lz^{\text{\rm DN}}=\lz^{\text{\rm DD}}=1$
with $g(x)=x$, $\big(\kz^{\text{\rm DN}}\big)^{-1}=\big(\kz^{\text{\rm DD}}\big)^{-1}\approx 2.1$.
For $\lz^{\text{\rm NN}}$, the symmetry in the whole line is lost. We have $\lz^{\text{\rm NN}}=2$
with $g(x)=-1+2 x^2$, $\big(\kz^{\text{\rm NN}}\big)^{-1}\approx 4.367$
which is achieved at $(x, y)\approx (0.316, 1.185)$.
Note that $\lim_{x\to\infty} g'(x)=\infty$  but $\lim_{x\to\infty}\big(e^C g'\big)(x)=0$.

To study ${\underline{\bar\kz}}^{\text{\rm NN}}$, recall that we can reduce
the NN-case to the DD-one by an exchange of $\mu$ and $\nu$.
By Proposition 3.1, we have $\big({\bar\kz}^{\text{\rm NN}}\big)^{-1}\approx 2.6$.
By Corollary {\ref{t033}} and
(\ref{3-36}), we obtain
$\big({\underline\kz}^{\text{\rm NN}}\big)^{-1}\approx 1.83$
with $(x^*, y^*, \uz^*)\approx (0.6405, 0.938, 0.721194)$.
For the last conclusion, we use a direct search star\-ting from
$(x, y)\approx (0.316, 1.185)$ which leads to $\kz^{\text{\rm NN}}$
in the last paragraph. The ratio becomes $2.6/1.83\approx 1.42<4$.
We mention that similar estimates can also be obtained by using a
different approximating procedure in parallel with \rf{cmf10}{Theorem 6.3}. Refer to \rf{cmf01}{Footnotes 12 and 14}.}\dexmp

The following examples are often illustrated in the textbooks on ordinary differential
equations, see for instance Hartman (1982), \S 11.1.

\xmp\lb{t042} {\rm The equation
$$u''+\sz^2 u=0\qd(\sz\ne 0)$$
has the general solution
$$u=c_1\cos(\sz x)+c_2 \sin(\sz x).$$
From this, it should be clear that
for the operator $L= \d^2/\d x^2$ with finite state space $(\az, \bz)$,
we have}
$$\aligned
\lz^{\text{\rm DD}}&= \bigg(\frac{\pi}{\bz-\az}\bigg)^2,\qqd
g(x)=\sin\bigg(\frac{\pi(x-\az)}{\bz-\az}\bigg);\\
\lz^{\text{\rm NN}}&= \bigg(\frac{\pi}{\bz-\az}\bigg)^2,\qqd
g(x)=\cos\bigg(\frac{\pi(x-\az)}{\bz-\az}\bigg);\\
\lz^{\text{\rm DN}}&= \bigg(\frac{\pi}{2(\bz-\az)}\bigg)^2,\qqd
g(x)=\sin\bigg(\frac{\pi(x-\az)}{2(\bz-\az)}\bigg);\\
\lz^{\text{\rm ND}}&= \bigg(\frac{\pi}{2(\bz-\az)}\bigg)^2,\qqd
g(x)=\cos\bigg(\frac{\pi(x-\az)}{2(\bz-\az)}\bigg).
\endaligned$$
{\rm The corresponding estimates are as follows.
$$\big(\kz^{\text{\rm DD}}\big)^{-1}
=\big(\kz^{\text{\rm NN}}\big)^{-1}
=\bigg(\frac{4}{\bz-\az}\bigg)^2,\qqd
\big(\kz^{\text{\rm DN}}\big)^{-1}
=\big(\kz^{\text{\rm ND}}\big)^{-1}
=\bigg(\frac{2}{\bz-\az}\bigg)^2.$$

Note that by symmetry, the DD- and NN-cases can be split at $\uz=(\az+\bz)/2$ into
the DN- and ND-cases. One can then approach $\lz^{\text{\rm DD}}$ and $\lz^{\text{\rm NN}}$ by using the known approximating method for $\lz^{\text{\rm DN}}$ and $\lz^{\text{\rm ND}}$ (cf. \rf{cmf01}{Theorem 1.2}). However, as an illustration of
Theorem {\ref{t032}} and Corollary {\ref{t033}}, we now compute ${\bar\kz}^{\text{\rm DD}}$ and ${\underline\kz}^{\text{\rm DD}}$.

Consider first the simpler interval $(\az, \bz)=(0, 1)$. Since $\mu=\nu=\d x$,
by symmetry, one may choose $y=1-x$. Then $x<1/2$ and
$$\aligned
\big({\bar\kz}^{\text{\rm DD}}\big)^{-1}\!\!&=
\!\inf_{x\in (0, 1/2)}\frac{2}{x}\bigg[1\!-\! 2 x
+ x^{-2}\int_0^x \! z^2\d z
+x^{-2}\int_{1-x}^1 \!(1-z)^2\d z
 \bigg]^{-1}\\
&=\inf_{x\in (0, 1/2)}\frac{6}{3 x -(2 x)^2}\\
&= \frac{32}{3}  \qd (\text{with }x=3/8).
\endaligned$$

To compute ${\underline\kz}^{\text{\rm DD}}$, set again $y=1-x$
with $x\in (0, 1/2)$. Then, the test function $f^{x, y}$ becomes
$$f^x(s)=
\begin{cases}
\sqrt{s \wedge x} &\qd s\le 1-x\\
\sqrt{1-s} &\qd s\in (1-x, 1).
\end{cases}
$$
By symmetry again, we have $\uz=1/2$. Fix $x\in (0, 1/2)$.
For convenience, we express $f^x$ as $(f_1, f_2)$: $f_1(s)=\sqrt{s}$ for $s\in [0, x]$ and
$f_2(s)=\sqrt{x}$ for $s\in [x, 1/2]$. Then by  (\ref{3-29})
with $\nu_-(s)=s$, we have $h^-=\big(h_1^-(z), h_2^-(z)\big)$:
$$\aligned
h_1^-(z)&=\int_0^z f_1(s) s \d s +z \bigg[\int_z^x f_1 +\int_x^{1/2} f_2 \bigg],\qd & z\in [0, x]\\
h_2^-(z)&=\bigg[\int_0^x f_1(s) s \d s + \int_x^z f_2(s) s\d s\bigg] +z\int_z^{1/2} f_2,\qd & z\in [x, 1/2].
\endaligned$$
Hence by (\ref{3-32}), we have
$$ I\!I^-(f^x)(z)=\frac{h^-(z)}{f^x(z)}=
\begin{cases}
\Big(-\frac 1 3 x^{3/2}+\frac 1 2 x^{1/2}\Big)\sqrt{z}-\frac{4}{15} z^{2}, &\qqd  z\in [0, x],\\
\frac{1}{10}\big(5 z (1-z)-x^2\big), & \qqd z\in [x, 1/2].
\end{cases}$$
Define
$$H(x)= -\frac 1 3 x^{3/2}+\frac 1 2 x^{1/2}\qqd\text{and}\qqd \gz (z)=H(x)\sqrt{z} - \frac{4}{15} z^{2}.$$
Then
$$\gz'(z)=\frac{H(x)}{2\sqrt{z}}-\frac{8}{15}z,
\qqd \gz''(z)=-\frac{H(x)}{4 z^{3/2}}- \frac{8}{15}<0.$$
Hence $\gz$ achieves its maximum at
$$z^*(x)=\bigg(\frac{15}{16} H(x)\bigg)^{2/3}.$$
Furthermore,
$$\gz(z^*(x))=H(x) \bigg(\frac{15}{16} H(x)\bigg)^{1/3}-\frac{4}{15}
\bigg(\frac{15}{16} H(x)\bigg)^{4/3}=\frac{3}{8}
\bigg(\frac{15}{2}\bigg)^{1/3} H(x)^{4/3}. $$
Note that $z^*(x)\le x$ iff $x\ge 5/14$. Besides, on the subinterval
$[x, 1/2]$, $h^-(z)/f^x(z)$ has maximum $1/8- x^2/10$ by (\ref{3-34}).
Solving the equation
$$\frac{3}{8}
\bigg(\frac{15}{2}\bigg)^{1/3} H(x)^{4/3}=\frac 1 8 - \frac{1}{10}
x^2,\qqd x\in (5/14, 1/2),$$
we obtain $x^*\approx 0.436273$ and then
$$\inf_{x\in (5/14, 1/2)}\,\sup_{z\le 1/2} \frac{h^-(z)}{f^x(z)}
=\gz(z^*(x^*))\approx 0.105967.$$
From these facts and (\ref{3-36}), we conclude that
$$\big({\underline \kz}^{\text{\rm DD}}\big)^{-1}
\approx 1/ 0.105967\approx 9.43693.$$
By the way, we mention that a similar but simpler study shows that
$$\inf_{x\in (0, 5/14)}\,\sup_{z\le 1/2} \frac{h^-(z)}{f^x(z)}
=\frac 1 8.$$
This shows that to get a less sharp lower bound $1/8$, the computation
becomes much simpler. It needs to study the extremal case that
$x=0$ only; the corresponding test function becomes
$f^x\equiv 1$.
Return to the original interval $(\az, \bz)$, by Proposition \ref{t031}{}
and Corollary \ref{t033}, we obtain
$$\frac{8}{(\bz-\az)^2}<\frac{9.4369}{(\bz-\az)^2}<
{\lz}^{\text{\rm DD}}= \bigg(\frac{\pi}{\bz-\az}\bigg)^2
\le \frac{32}{3(\bz-\az)^2}=\frac 2 3 \bigg(\frac{4}{\bz-\az}\bigg)^2.$$
The ratio becomes $\frac{32}{3}\big/ 9.4369\approx 1.13$.
The same assertion holds if ${\lz}^{\text{\rm DD}}$ is replaced by ${\lz}^{\text{\rm NN}}$
because of the symmetry.

 It is a good chance to discuss the approximating procedure remarked after
Corollary {\ref{t033}}. Here we consider the lower estimate only.
Replacing $f^x=(f_1, f_2)$
by $(h_1^-, h_2^-)$, one produces a new $(h_1^-, h_2^-)$ and then a new $I\!I^-(f)$ which provides a new lower bound. By using this procedure twice with fixed $\uz=1/2$ and $x=x^*\approx 0.436273$, we obtain successively the following lower bounds:
$$\frac{9.80392}{(\bz-\az)^2}, \qqd \frac{9.86193}{(\bz-\az)^2}.$$
 Clearly, they are quite close to the exact value of ${\lz}^{\text{\rm DD}}$ and ${\lz}^{\text{\rm NN}}$:
$$\frac{\pi^2}{(\bz-\az)^2}\approx \frac{9.8696}{(\bz-\az)^2}.$$}\dexmp

\xmp{} {\rm By a substitute $u=z e^{-b x/2}$, the
equation
$$u''+b u' +\gz u=0\qqd(b, \gz \text{ are real constants})$$
is reduced to
$$z''+\sz^2 z=0 \qqd \big(\sz^2= \gz-b^2/4\big).$$
From the last example, it follows that the equation has general solutions
$$u=
\begin{cases} e^{-b x/2} (c_1+c_2 x) \; &\text{if } \gz=b^2/4\\
c_1 e^{\xz_1 x}+c_2 e^{\xz_2 x}\; &\text{if } \gz<b^2/4\\
e^{-b x/2} \Big(c_1 \cos\big(x \sqrt{\gz-b^2/4}\big) +c_2 \sin\big(x
\sqrt{\gz-b^2/4}\,\big)\Big) &\text{if $\gz>b^2/4$},
\end{cases}$$
where $\xz_1, \xz_2$ are solution to the equation
$$\xz^2+b\,\xz +\gz=0.$$
Thus, for the operator $L=\d^2/\d x^2 + b\,\d/\d x$\,($b$ is a constant) with
state space $(0, \infty)$, we have the following principal
eigenfunctions
\begin{itemize}
\setlength{\itemsep}{-0.8ex}
\item $g(x)=(2/b+x) e^{-b x/2}$ and $g(x)=x e^{-b x/2}$ in ND- and DD-cases,
respectively, when $b>0$;
\item $g(x)=x e^{-b x/2}$ and $g(x)=(1+b x/2)e^{-b x/2}$ in DN- and NN-cases,
respectively, when $b<0$.
\end{itemize}
In each of these cases, we have the principal eigenvalue $\lz^{\#}=b^2/4 $ and $\big(\kz^{\#}\big)^{-1}=b^2$. Moreover,
$\big(\bar\kz^{\text{\rm DD}}\big)^{-1}$\!, $\big(\bar\kz^{\text{\rm NN}}\big)^{-1}=b^2/2$.
Clearly, the lower estimate $\big(\kz^{\#}\big)^{-1}/4$ is sharp
in all cases.}
\dexmp

\xmp [Cauchy--Euler equation]{\rm Consider the operator
$$L=x^2 \frac{\d^2}{\d x^2} + b x \frac{\d}{\d x},$$
where $b$ is a constant. By a change of
variable $x=e^{y}$, the equation
$$x^2 u''+b x u'+\gz u=0\qqd\text{($b$, $\gz$ are constants)}$$
is reduced to the last example:
$$\frac{\d^2 u}{\d y^2}+(b-1)\frac{\d u}{\d y}+\gz u =0.$$
Hence the original equation has general solutions
$$u\!=\!
\begin{cases}
{x^{(1-b)/2}}\, (c_1+ c_2 \log x)&\qqd\qqd\text{if $\gz=(1-b)^2/4$}\\
c_1 x^{\xz_1}+c_2 x^{\xz_2} &\qqd\qqd \text{if $\gz<(1-b)^2/4$}\\
{x^{(1-b)/2}} \Big(c_1 \cos\big(\!\sqrt{\gz-(1-b)^2/4}\,\log x\big)\! +\! c_2
&\!\!\!\!\sin \big(\!\sqrt{\gz-(1-b)^2/4}\,\log x\big)\!\Big)\\
&\qqd\qqd\text{if $\gz>(1-b)^2/4$,}
\end{cases}$$
where $\xz_1, \xz_2$ are solution to the equation $\xz^2+
(b-1)\xz+\gz=0:$
$$\xz_1, \xz_2=(1-b)/2 \pm \sqrt{(1-b)^2/4-\gz} .$$
Here we have used Euler's formula:
$$x^{i\sqrt{\xz}}=e^{i\sqrt{\xz}\,\log x}=\cos\big(\sqrt{\xz}\log x\big)
+ i \sin\big(\sqrt{\xz}\log x\big).$$
In particular,

(1) when $b=0$, we have solutions
$$u=
\begin{cases}
\sqrt{x}\, (c_1+ c_2 \log x)&\qd\text{if $\gz=1/4$}\\
c_1 x^{\xz_1}+c_2
 x^{\xz_2}&\qd\text{if $\gz<1/4$}\\
\sqrt{x}\Big(c_1 \cos\big(\sqrt{\gz-1/4}\,\log x\big)
\!+\!c_2
\sin\big(\sqrt{\gz-1/4}\,\log x\big)\Big)&\qd\text{if $\gz>1/4$}.
\end{cases}$$
Now, corresponding to $\gz=1/4$, we have
$$\lz^{\text{\rm DN}}=\frac 1 4, \qqd g(x)={\begin{cases}
\sqrt{x} &\text{ if the state space is }(0, \infty)\\
\sqrt{x}\,\log \sqrt{x}&\text{ if the state space is }(1, \infty).
\end{cases}}$$
The first case is the original Hardy's inequality.
Corresponding to $\gz=1/4$ again but for state space $(1, \infty)$, we have
$$\lz^{\text{\rm NN}}=\frac 1 4, \qqd g(x)=\sqrt{x}\,\big(\log \sqrt{x}-1\big).$$
Here $\lim_{x\to\infty}\big(e^C g'\big)(x)=\lim_{x\to\infty} g'(x)=0$.
We have $\big(\kz^{\text{\rm DN}}\big)^{-1}$,
$\big(\kz^{\text{\rm NN}}\big)^{-1}=1$,
$\big(\bar\kz^{\text{\rm DN}}\big)^{-1}$,
$\big(\bar\kz^{\text{\rm NN}}\big)^{-1}=1/2$,
respectively.
The lower estimate $\big(\kz^{\#}\big)^{-1}/4$ is sharp in each case.
The DN-case is actually a special one of the last example.

(2) When $b=1$, for finite state space $(1, N)$ with Dirichlet boundaries, we have
$$\aligned
\lz_n&= \bigg(\frac{n\pi}{\log N}\bigg)^2,\qqd
g(x)=\sin\bigg(\frac{n\pi}{\log N}\log x\bigg),\qqd n\ge 1.
\endaligned$$
In particular,
$$\lz^{\text{\rm DD}}=\bigg(\frac{\pi}{\log N}\bigg)^2,\qqd
g(x)=\sin\bigg(\frac{\pi}{\log N}\log x\bigg).$$
Next, for Neumann boundaries, we have
$$\lz^{\text{\rm NN}}= \bigg(\frac{\pi}{\log N}\bigg)^2,\qqd
g(x)=\cos\bigg(\frac{\pi}{\log N}\log x\bigg).$$
In both cases, we have
$\big(\kz^{\text{\rm DD}}\big)^{-1}$,
$\big(\kz^{\text{\rm NN}}\big)^{-1} =\big({4}/{\log N}\big)^2$.
Besides, we have
$$\aligned
\lz^{\text{\rm DN}}&= \bigg(\frac{\pi}{2\log N}\bigg)^2,\qqd
g(x)=\sin\bigg(\frac{\pi}{2\log N}\log x\bigg);\\
\lz^{\text{\rm ND}}&= \bigg(\frac{\pi}{2\log N}\bigg)^2,\qqd
g(x)=\cos\bigg(\frac{\pi}{2\log N}\log x\bigg).\endaligned$$
In these cases, we have $\big(\kz^{\text{\rm DN}}\big)^{-1}$,
$\big(\kz^{\text{\rm ND}}\big)^{-1} =\big({2}/{\log N}\big)^2$.
Note that the present case can be reduced to Example \ref{t042} under
the change of variable $x=e^y$, the results here can be obtained
from Example \ref{t042} replacing $(\az -\bz)^2$ by $\log^2N$.
In view of this, we also have
$$\big(\bar\kz^{\text{\rm DD}}\big)^{-1}=
\big(\bar\kz^{\text{\rm NN}}\big)^{-1}=\frac{32}{3\log^2 N},\qqd
\big(\underline\kz^{\text{\rm DD}}\big)^{-1}=
\big(\underline\kz^{\text{\rm NN}}\big)^{-1}
\approx\frac{9.4369}{\log^2 N}.$$}
\dexmp

\section{Appendix}

The next result is a generalization of \rf{cmf10}{Proposition 1.2}.

\prp {\cms Let $P_t(x, \cdot)$ be symmetric and have density
$p_t(x, y)$ with respect to $\mu$. Suppose that the diagonal elements
$p_{s}(\cdot, \cdot)\in L_{\text{\rm loc}}^{1/2}(\mu)$
for some $s>0$ and a set $\scr K$ of bounded functions with compact support is
dense in $L^2(\mu)$. Then $\lambda_0 = \varepsilon_{\max}$.
}
\deprp

\prf The proof is similar to the ergodic case (cf.
\rf{cmf05a}{Section 8.3} and \rf{cmf10}{proof of Theorem 7.4}), and is
included here for completeness.

(a) Certainly, the inner product and norm here are taken with respect to $\mu$. First, we have
$$\aligned
P_t(x, K)&= P_s P_{t-s} \mathbbold{1}_K (x)\\
&=\int \mu(\d y) \frac{\d P_s(x, \cdot)}{\d \mu}(y) P_{t-s} \mathbbold{1}_K (y)\qd\text{(since $P_s\ll \mu$)}\\
&=\mu\bigg(\frac{\d P_s(x, \cdot)}{\d \mu} P_{t-s} \mathbbold{1}_K\bigg)\\
&=\mu\bigg(\mathbbold{1}_K P_{t-s}\frac{\d P_s(x, \cdot)}{\d \mu} \bigg)\qd\text{(by symmetry of $P_t$)}\\
&\le \sqrt{\mu(K)}\,\bigg\|P_{t-s}\frac{\d P_s(x, \cdot)}{\d \mu}\bigg\|\qd\text{(by Cauchy-Schwarz inequality)}\\
&\le \sqrt{\mu(K)}\,\bigg\|\frac{\d P_s(x, \cdot)}{\d \mu}\bigg\|\, e^{-\lz_0 (t-s)} \qd\text{(by $L^2$-exponential convergence)}\\
&=\Big(\sqrt{\mu(K)\, p_{2s}(x, x)}\,e^{\lz_0 s}\Big) e^{-\lz_0 t}\qd\text{(by \rf{cmf05a}{(8.3)})}.
\endaligned$$
By assumption, the coefficient on the right-hand side is locally $\mu$-integrable.
This proves that $\vz_{\max}\ge \lz_0$.

(b) Next, for each $f\in{\scr K}$ with $\|f\|=1$, we have
\[\aligned
\|P_t f\|^2&= (f, P_{2t} f)\qd\text{(by symmetry of $P_t$)}\\
&\le \|f\|_{\infty}\int_{\text{\rm supp}\,(f)} \mu(\d x) P_{2t} |f|(x)\\
&\le \|f\|_{\infty}^2 \int_{\text{\rm supp}\,(f)} \mu (\d x) P_{2t} (x, \text{\rm supp}\,(f))\\
&\le \|f\|_{\infty}^2 \int_{\text{\rm supp}\,(f)} \mu (\d x) c(x, \text{\rm supp}\,(f))
e^{-2\vz_{\max} t}\\
&=: C_f e^{-2\vz_{\max} t}.
\endaligned\]
The technique used here goes back to Hwang et al. (2005).

(c) The constant $C_f$ in the last line can be removed.
Following Lemma 2.2 in Wang (2002),
by the spectral representation theorem and the fact that $\|f\|=1$, we have
\[\aligned
\|P_t f\|^2&=\int_0^\infty e^{-2 \lambda t}\d (E_\lambda f, f)\\
&\ge \bigg[\int_0^\infty e^{-2 \lambda s}\d (E_\lambda f,
f)\bigg]^{t/s}\quad \text{(by Jensen's inequality)}\\
&=\|P_s f\|^{2t/s},\qquad \;t\ge s.\endaligned\]
Note that here the semigroup is allowed to be sub-Markovian.
Combining this with (b), we have
$\|P_s f\|^2\le C_f^{s/t}
e^{-2 \vz_{\max} s}$. Letting $t\to \infty$, we obtain
$$\|P_s f\|^2\le e^{-2\vz_{\max} s}, $$
first for all $f\in {\scr K}$ and then for all $f\in L^2(\mu)$
with $\|f\|=1$, because of the denseness of ${\scr K}$ in $L^2(\mu)$.
Therefore, $\lambda_0\ge \varepsilon_{\max}$.
Combining this with (a), we complete the proof.
\deprf
\medskip

The main result (Theorem \ref{t2}) of this paper is presented in the last section
(section 10) of the paper \ct{cmf10}, as an analog of birth-death processes. Paper \ct{cmf10}, as well as \ct{cmf08} for $\fz^4$-model,  is available on arXiv.org.

\medskip{\small

\nnd{\bf Acknowledgments}.
Research supported in part by the Creative Research Group Fund of the National
Natural Science Foundation of China (No. 10721091), by the ``985'' project from the Ministry of Education in China.
The author has been luckily invited by Professor Louis Chen three times
with financial support to visit Singapore. Deep appreciation is given to him for
his continuous encouragement and friendship in the past 30 years.
Sections \ref{s0}--\ref{s3} of the paper are based on the talks presented in
``Workshop on Stochastic Differential Equations and Applications''
(December, 2009, Shanghai),
``Chinese-German Meeting on Stochastic Analysis and Related Fields''
(May, 2010, Beijing), and ``From Markov Processes to Brownian Motion and Beyond ---
An International Conference in Memory of Kai-Lai Chung'' (June, 2010, Beijing).
The author acknowledges the organizers of the conferences: Professors
Xue-Rong Mao; Zhi-Ming Ma and Michael R\"okner; and the Organization Committee headed by
Zhi-Ming Ma (Elton P. Hsu and Dayue Chen, in particular),
for their kind invitation and financial support.}


\nnd {\small School of Mathematical Sciences, Beijing Normal University,
Laboratory of Mathematics and Complex Systems (Beijing Normal University),
Ministry of Education, Beijing 100875,
    The People's Republic of China.\newline E-mail: mfchen@bnu.edu.cn\newline Home page:
    http://math.bnu.edu.cn/\~{}chenmf/main$\_$eng.htm
}
\end{document}